# J. K. Ghosh's contribution to statistics: A brief outline

## Bertrand Clarke[1] and Subhashis Ghosal[2]

*University of British Columbia and North Carolina State University*

**Abstract:** Professor Jayanta Kumar Ghosh has contributed massively to various areas of Statistics over the last five decades. Here, we survey some of his most important contributions. In roughly chronological order, we discuss his major results in the areas of sequential analysis, foundations, asymptotics, and Bayesian inference. It is seen that he progressed from thinking about data points, to thinking about data summarization, to the limiting cases of data summarization in as they relate to parameter estimation, and then to more general aspects of modeling including prior and model selection.

## Contents

[1]Department of Statistics, University of British Columbia, 6356 Agricultural Road, Vancouver, BC, V6T1Z2, Canada, e-mail: riffraff@stat.ubc.ca
[2]Department of Statistics, North Carolina State University, 12 Patterson Hall, 2501 Founders Drive Raleigh, NC 27695, USA, e-mail: ghoshal@stat.ncsu.edu

## 1. Introduction

Professor Jayanta Kumar Ghosh, or J. K. Ghosh, as he is commonly known, has been a prominent contributor to the discipline of statistics for five decades. The spectrum of his contributions encompasses sequential analysis, the foundations of statistics, finite populations, Edgeworth expansions, second order efficiency, Bartlett corrections, noninformative, and especially matching, priors, semiparametric inference, posterior limit theorems, Bayesian nonparametrics, model selection Bayesian hypothesis testing and high dimensional data analysis, as well as some applied work in reliability theory, statistical quality control, modeling hydrocarbon discoveries, geological mapping and DNA fingerprinting. By itself, covering such diverse topics in depth is a major career achievement. He has authored over 130 publications including three monographs and several edited volumes. His books, one entitled *Higher Order Asymptotics* and published as an IMS monograph and another entitled *Bayesian Nonparametrics*, co-authored by R. V. Ramamoorthi and published by Springer-Verlag, continue to hold respected positions for researchers in these areas. His recently published third book [34] is a fine graduate text on Bayesian inference. In addition, his service to the profession, especially as the editor of Sankhyā, has been invaluable.

The variety of his work notwithstanding, asymptotics has been central to his thinking across a wide range of problems. Accordingly, in what follows, we outline some of his work, in roughly chronological order, focussing on those contributions which are intimately connected to asymptotics. In the course of reviewing his work, we try to characterize the progression of thinking that naturally connects the topics that J. K. Ghosh has done so much to develop.

## 2. Sequential analysis

J. K. Ghosh started his research career in Sequential Analysis in the early sixties as a Graduate Student in the Department of Statistics at Calcutta University. Wald had recently introduced his *sequential probability ratio test* (SPRT), but its properties were not well understood in the composite case. This was the first topic to which Ghosh turned his attention. Through his work, many of the properties of SPRT and related procedures were established and better understood. For instance, in the testing context, double minimaxity essentially means simultaneous minimization of average type I and II error probabilities. In his first published work [26], Ghosh clarified a result of Wald on the double minimaxity of the SPRT for normal two-sided alternative hypothesis (with unknown scale) separated from the null by $\delta$.

It is well-known that the power function is monotonic in many common families for fixed sample sizes. Ghosh established an analog of this result in [27], namely that the operating characteristic function of the (generalized) SPRT continues to be monotonic. Also in the sequential context, [28] considered the admissibility of sequential tests based on a simple identity which later became known as the Ghosh–Pratt identity. Ghosh compared the SPRT not just with the class of all tests with finite expected sample size but also within other classes, for instance, the class which requires at least one observation or which requires no more than a predetermined number of observations to reach a conclusion.

Following this, Ghosh continued to elucidate more properties of the SPRT, and its variants, which could be seen as analogs of the corresponding properties Neyman–Pearson or Bayes tests for fixed sample size. In [29], he proved that for exponential families, truncated or untruncated Bayesian sequential decision rules' terminal decisions describe regions in terms of sufficient statistics, and also showed that for testing problems, truncated generalized SPRT's form a complete class.

About two decades later, Ghosh returned to sequential problems, along with various co-authors. In [33], he studied an invariant SPRT to identify two normal populations with equal variance and obtained bounds for error probabilities. Most recently, similar bounds for an invariant SPRT with respect to an improper prior have also been obtained in [50].

Two-stage procedures are closely related to sequential procedures. Recall Stein's famous problem of finding a bounded length confidence interval for the normal mean with unknown variance. Stein proposed a two-stage procedure for doing this: In the first stage, the sample variance determines how many samples are to be taken in the second stage. An obvious shortcoming of the procedure is that the second stage sample variance is not used in the construction of the interval. So, it is natural to ask whether one can improve Stein's procedure by using the second stage sample variance. Surprisingly, it is impossible to better Stein's procedure, as shown in [38].

However, the procedure can be improved in a different, and perhaps more appropriate sense. The confidence coefficient does not in general properly reflect the true sense of confidence about a parameter after observing data. For instance, if two observations are obtained from a $U(\theta, \theta+1)$ family, then the assessment of $\theta$ is very precise when the two observations differ in magnitude by nearly 1, while the assessment is much less precise if the two observations are close to each other. This means that classical confidence intervals fail to indicate the true difference in the level of confidence after observing the sample.

Motivated by this, Kiefer suggested letting the confidence coefficient depend on the data. After all, in reality, for a given random interval $I$, we often want to predict the indicator function $1\{\theta \in I\}$. Since this object is unknown, it is traditionally estimated by a constant, the best constant being the expectation $P_\theta(\theta \in I)$, which becomes fixed (or asymptotically fixed) for many classical intervals. However, from a prediction theory point of view, it makes more sense to let the predictor of $1\{\theta \in I\}$ depend on the observed data. The predictor considered in this way is called the random confidence coefficient associated with the confidence interval $I$. It is shown in [39] that the second stage sample variance can be used to boost the random confidence coefficient of a bounded length confidence interval.

## 3. Foundations of statistics

From the examination of individual data points as they relate to the testing problem, Ghosh shifted his attention to data summarization, focussing on the relationship between sufficiency and invariance. Sufficiency isolates features of the collection of observations from those of the individual ones which are independent of the features of the collection. Invariance, on the other hand, summarizes data by imposing symmetry constraints. In practice, both sufficiency and invariance restrictions are applied, but their order of application is an issue of interest.

Consider a statistical model $(\mathfrak{X}, \mathscr{A}, \mathcal{P})$ where a group of transformations $G$ is acting on the sample space and attention is limited to invariant procedures. To find a sufficiency reduction, one needs to find a sufficient sub-$\sigma$-field $\mathscr{C}$ of the invariant $\sigma$-field $\mathscr{I}$. However, in practice, it is typically easier to invoke invariance on the data after it has been reduced by sufficiency. Let $\mathscr{S}$ be a sufficient $\sigma$-field. To justify the application of invariance restriction after a sufficiency reduction, it is enough to establish that $\mathscr{S} \cap \mathscr{I}$ is sufficient for $\mathscr{I}$. This problem was addressed by W. J. Hall, R. A. Wijsman and J. K. Ghosh, independently and roughly simultaneously. Once they realized they had compatible results, they published a combined paper [65]. Their main result can be described briefly as follows. A statistic $T$ is called almost invariant if, for every $g \in G$, $T(x) = T(gx)$ a.s. Under conditions that imply that every almost invariant set is equivalent, up to null sets, to an invariant set, it follows that $\mathscr{S}$ and $\mathscr{I}$ are conditionally independent given $\mathscr{S} \cap \mathscr{I}$, and hence $\mathscr{S} \cap \mathscr{I}$ is sufficient for $\mathscr{I}$.

Another notion which relates two sequences of $\sigma$-fields in sequential experiments is that of transitivity, introduced by Bahadur. Two sequences of $\sigma$-fields $\mathscr{B}_n \subset \mathscr{A}_n$ are said to be transitive if for every $\mathscr{B}_{n+1}$-measurable function $f$, $\mathrm{E}(f|\mathscr{A}_n)$ is $\mathscr{B}_n$ measurable. In the usual sequential setting up, $\mathscr{S}_n \cap \mathscr{I}_n$ is transitive for $\mathscr{I}_n$, where the extra index $n$ indicates the sample size. Several implications of this result were discussed in [65].

In many application areas, sample surveys for instance, discrete models arise, where the probability is concentrated on a countable set but the models do not have common support, i.e., the support set is different for different parameter values. Clearly, such a family is not dominated and the Halmos–Savage theorem on sufficiency does not hold. Nevertheless, as shown in [2], minimal sufficient $\sigma$-fields exist and the Neyman factorization theorem holds good. These results were extended for pairwise sufficient $\sigma$-fields and condition for existence of minimal pairwise sufficient $\sigma$-field was found in [37].

Another basic question is whether a fixed-dimensional sufficient statistic independent of sample size actually exists. In exponential families, it is well-known that fixed-dimensional sufficient statistics exist. Outside of exponential families, however, sufficient statistics are hard to find. Some distinguished nonregular cases like $U(-\theta, \theta)$ provide additional examples. In [54], it is shown that if the support $(a(\theta), b(\theta))$ is shrinking or expanding, as in the support of $U(0, \theta)$ for example, then the density must be of the form $g(\theta)h(x)$ to have a real-valued sufficient statistic. If $a(\theta)$ and $b(\theta)$ are both increasing, or both decreasing, as in $U(\theta, \theta+1)$, then an $\mathbb{R}^2$-valued sufficient statistic can exist only in special cases.

## 4. Asymptotics

The asymptotic point of view undergirded Ghosh's thinking, even in problems that were not primarily focussed on asymptotic properties. In a sense, much of his work on sequential analysis, Bayesian analysis and Bayesian nonparametrics are also, at least implicitly, work on asymptotics. In fact, many of the most important asymptotic ideas, such as higher order asymptotics and Edgeworth expansions, were pioneered by him. Moreover, in terms of how his thinking progressed, asymptotics can be regarded as the next natural conceptual step after thinking about data points in sequential analysis, and sufficiency or invariance as a data summarization technique. That is, once we have gathered and summarized our data, we want to see where it seems to be leading us.

Ghosh's work on asymptotics can be broadly grouped into seven categories. He worked on the Bahadur–Ghosh–Kiefer representation for a quantile. He made foundational contributions to establishing the existence of Edgeworth expansions. In higher order asymptotics, Ghosh examined second order efficiency, Bartlett correction and contributed to our understanding of how Wald, Rao and likelihood ratio tests compare. Then he turned his attention to Bahadur efficiency and the vexing Neyman–Scott problem.

### 4.1. Bahadur–Ghosh–Kiefer representation

Bahadur represented a sample quantile approximately as an average of i.i.d. random variables. To get this representation, Bahadur assumed the existence of two derivatives of the c.d.f.; the second derivative is bounded and the first derivative is positive on a neighborhood of the $p$-th population quantile $\xi_p$. Then Bahadur showed that the error in the representation is $O(n^{-3/4}(\log n)^{3/4})$. The order of error in Bahadur's representation is nearly sharp, cf. the exact order $n^{-3/4}(\log n)^{1/2}(\log \log n)^{1/4}$ obtained by Kiefer.

One of the reasons this is important is that the order of error is small enough to obtain asymptotic normality for the sample quantiles. However, assuming the existence of two derivatives is somewhat strong. For instance, it rules out the location family from the double exponential density. On the other hand, for most statistical purposes, where only the asymptotic distribution is important, having an error term of order $o_p(n^{-1/2})$ is enough. Therefore it is of interest to weaken Bahadur's assumptions at the expense of weakening the conclusion to $o_p(n^{-1/2})$. This is possible, even for a variable point $p_n$ depending on $n$, as shown in [30], using only the assumption of positive first derivative at $\xi_p$.

Actually, the idea of representing a quantile approximately as an average of i.i.d. observations occurred to Ghosh independently in the mid sixties at the same time as Bahadur was working on the problem. Ghosh was looking at the problem in the more general multivariate multisample framework in connection with asymptotic normality of multivariate rank tests. He did not record his proof then since it did not extend to the multivariate setup at that time.

### 4.2. Edgeworth expansions

Edgeworth expansions are natural refinements of asymptotic normality results in that they give error terms of asymptotically smaller order by including more terms in addition to the leading normal term. However, for a long time, Edgeworth expansions were only heuristically justified. In the pioneering paper [7], it is shown that under conditions of finiteness of certain moments and a condition on characteristic function known as Cramér's condition is the literature, the $r$-th order Edgeworth expansion of a smooth function of sample averages admits error $O(n^{-(r+1)/2})$. In particular, it follows that for the sample average, finiteness of the $2r$-th moments is required to justify an Edgeworth expansion of order $r$. The follow-up paper [8] relaxes some moment conditions. A thorough and lucid treatment of Edgeworth expansions and higher order asymptotics is given in Ghosh's IMS monograph [32].

Another angle on Edgeworth expansions comes from Fisher consistency. Consider an exponential family with density proportional to $\exp[\sum_{j=1}^{k} w_j(\theta)t_j(x)]$. In this context, an estimator $T_n$, which is a function of the $k$-dimensional sufficient statistic $(\sum_{i=1}^{n} t_1(X_i), \ldots, \sum_{i=1}^{n} t_k(X_i))$, is Fisher consistent if $T_n(w(\theta)) = \theta$. Assuming sufficient smoothness conditions and linear independence of the component functions $w_1(\theta), \ldots, w_k(\theta)$, a Fisher consistent estimator can be written as a smooth function of sample averages, and hence has an Edgeworth expansion. In [62], this Edgeworth expansion is compared with that of the MLE, which is another Fisher consistent estimator. Interestingly, for any bowl shaped loss function, the MLE has better second order risk properties than any other Fisher consistent estimator. Consequently, this gives a way to discriminate among estimators which are first order asymptotically equivalent. This property is called second order efficiency and will be discussed in the next subsection.

Edgeworth-type expansions need not be restricted to asymptotically normal estimators. Other limiting distributions can appear naturally. Recall that log-likelihood ratio type statistics are among the most common statistics converging to non-normal limits such as a chi-square distribution. For locally quadratic functions of sample averages, such as the log-likelihood ratio, asymptotic expansions have been obtained in [13]. They have a leading chi-square term. Subsequent terms appear as coefficients of powers of $n^{-1/2}$ and are finite linear combinations of chi-square distributions of degrees of freedoms $p$, $p+2$, $p+4$, etc., where $p$ is the degree of freedom of the leading term. Similar expansions hold even under contiguous alternatives with non-central chi-squares replacing the chi-square leading term as shown in [14]. The subsequent terms are finite linear combinations of non-central chi-square distributions with degrees of freedoms $p$, $p+2$, $p+4$, and so forth.

### 4.3. Second order efficiency

Second order efficiency (also called third order efficiency by some authors) is the natural way to compare two asymptotically efficient estimators since they are first order equivalent. In particular, it was widely believed that the MLE, or some suitable variant of it, had, asymptotically, the smallest possible risk up to the second order. Ghosh, among others like Efron, Chibisov, Pfanzagl, Akahira and Takeuchi, made pioneering contributions towards rigorous justification of this assertion in [64]. His main result may be roughly described as follows: Let $T_n$ be an efficient estimator and consider a modification $T'_n = T_n + m(T_n)/n$. Then $T'_n$ can be beaten by $\hat{\theta}'_n = \hat{\theta}_n + g(\hat{\theta}_n)/n$, a modification of the MLE $\hat{\theta}_n$, where the function $g$ depends on $T_n$ and $m$. Here, by a better estimator, we mean that

$$\lim_{n \to \infty} n^2 [\mathrm{E}_\theta \{W(T'_n, \theta)\} - \mathrm{E}_\theta \{W(\hat{\theta}'_n, \theta)\}] \geq 0,$$

for all $\theta \in \Theta$, for a truncated squared error loss $W$. This paper also contains other impressive results such as Bhattacharya-type bounds, a Bayesian connection with second order efficiency and a notion of second order asymptotic sufficiency. Similar results about second order efficiency of the MLE for Pitman closeness and any bounded bowl shaped loss functions are given in [63].

In addition to second order efficiency, there is a notion of second order admissibility. An estimator is second order admissible if there is no estimator which has uniformly smaller second order risk with strict inequality for at least one point. In [59], for estimators of the form $\hat{\theta}_n + g(\hat{\theta}_n)/n$, a necessary and sufficient condition for second order admissibility under squared error loss is obtained.

These second order optimality properties of modified versions of the MLE raise the issue whether the MLE has optimality properties beyond the second order. A nice counterexample in [60], however, concludes negatively. On the other hand, questions on second order admissibility go beyond the MLE to any BAN estimator $\hat{\theta}_n$ modified to $\hat{\theta}_n + g(\hat{\theta}_n) + o_p(n^{-1})$. The condition for second order Pitman admissibility is obtained in [58], and its multiparameter version in [49].

Another natural question in the context of second order admissibility is the following. If two or more statistics are separately second order admissible, for two different components of a parameter with bias $o(n^{-1})$, then, is it true that they are jointly second order admissible? The question has a curious answer given in [16]. For two dimensions, they are jointly second order admissible, but for three or more dimensions, they are not jointly second order admissible. This result is reminiscent of Stein's phenomenon on ordinary admissibility with respect to the squared error loss for estimating the normal mean. Intuitively, asymptotically, all regular experiments are normal experiments and thus a phenomenon under normality continues to hold asymptotically under any regular model. The interesting part of the result is that the phenomenon shows up in the second order.

### *4.4. Bartlett correction*

Bartlett introduced a remarkable technique, which bears his name, to improve the chi-square approximation to the distribution of a log-likelihood ratio statistic. The idea is embarrassingly simple: rescale the chi-square distribution with the second order expansion of the mean of the statistic. It is surprising that such a simple strategy improves the approximation so much.

In the seminal paper [9], a variant on Wilks' theorem tuned to the goal of understanding the Bartlett correction was presented. Recall that Wilks' theorem is the statement that the log-likelihood ratio is asymptotically chi-square. However, the chi-square is the result of squaring normals. To see how this might apply to the log-likelihood ratio statistic, let $X_1, X_2, \ldots$ be i.i.d. observations from a parametric family governed by $\theta = (\theta^1, \ldots, \theta^p)$ and let $L(\theta)$ be the log likelihood. For $j = 1, \ldots, p$, let $\hat{\theta}_j$ be the MLE of $\theta$ under the null hypothesis $\theta^1 = \theta_0^1, \ldots, \theta^j = \theta_0^j$, and let $T_j = 2n(L(\hat{\theta}_{j-1}) - L(\hat{\theta}_j))^{1/2}\text{sign}(\hat{\theta}_{j-1} - \hat{\theta}_j)$, where $\hat{\theta}_0$ stands for the unrestricted MLE. Note that squaring $T_j$ gives the usual object in Wilks' theorem with limiting chi-square behavior. However, now, without the square, $(T_1, \ldots, T_p)$ is asymptotically normal with error $O_p(n^{-3/2})$ under the grand null hypothesis. This property of $T_j$ gives rise to the Bartlett correction in the multidimensional setting.

Another result developed in that paper is a Bayesian version of the Bartlett correction. This is a Bartlett correction to the posterior distribution, conditional on the data, obtained by letting the prior tend to the degenerate distribution at the true parameter value. The relation between the Bartlett correction and the Bayesian correction gives a deeper understanding of the Bartlett correction phenomenon and leads to a variety of generalizations.

Following this path, [41] studied the asymptotic equivalence of the frequentist and Bayesian Bartlett corrections for the likelihood ratio and the conditional likelihood ratio statistic (CLR) introduced by Cox and Reid. In particular, the conditions for equivalence are instrumental for giving a simple proof of the existence of the frequentist Bartlett correction for the CLR statistic. This was extended to the multivariate case in [40]. A variant on the likelihood ratio called the adjusted likelihood ratio (ALR) was introduced by McCullagh and Tibshirani. In [45], it was shown that the ALR statistic has behavior similar to that of the CLR statistic, in that it admits a Bartlett correction and its power under contiguous alternatives is equivalent to that of the CLR up to the order $o(n^{-1/2})$. In terms of average power, the agreement continues up to $o(n^{-1})$.

### *4.5. Comparison of the likelihood ratio, Wald's and Rao's statistics*

The problem of comparing the likelihood ratio (LR), Wald's and Rao's tests, with regard to power has received significant attention in the statistics and econometrics literature. It is well-known that, up to the first order of approximation and under contiguous alternatives, these three tests have the same local power as dictated by the noncentral chi-square distribution. Discrimination among them, therefore, calls for comparison via higher order power. However, while the LR test is locally unbiased up to a higher order of approximation, the same does not hold in general for the other two tests. From this perspective, to make them really comparable, Ghosh suggested considering locally unbiased versions of Wald's and Rao's tests. This work, done under his supervision, eventually led to an optimum property of Rao's test in terms of third order local power. A review of these developments is available in [31].

In addition, the power properties of the three tests as well as their Bartlett adjustability, when they are developed on the basis of a quasi-likelihood rather than a true density-based likelihood was discussed in [48].

### 4.6. Bahadur–Cochran deficiency

To compare the asymptotic performance of two tests, one may look at their Bahadur–Cochran relative efficiency, which is the limit, as $\delta \to 0$, of the ratio of the smallest integers which make their levels less than $\delta$. For many pairs of reasonable tests, the ratio turns out to be 1. To compare them at a finer level, it is sensible to look at their difference, which may be called the Bahadur–Cochran deficiency. The limit inferior (or superior) of the difference, reflecting the relative advantage of one test over other, was calculated in [12] for some common test statistics.

### 4.7. Neyman–Scott problem and semiparametric inference

Ghosh made notable contributions in the Neyman–Scott problem also. In the Neyman–Scott problem, a new observation $X_i$ is governed by a common parameter $\theta$ and an additional parameter $\xi_i$, depending on $i$, but only the parameter $\theta$ is of interest. The problem is notoriously difficult in that common estimators, such as the MLE, are usually inconsistent. For instance, if the $X_{ij}$ are independently normally distributed with mean $\mu_i$ and variance $\sigma^2$, $i = 1, \ldots, n$, $j = 1, \ldots, k$, then the MLE of $\sigma^2$, $(nk)^{-1} \sum_{i=1}^n \sum_{j=1}^k (X_{ij} - \bar{X}_i)^2$, where $\bar{X}_i = k^{-1} \sum_{j=1}^k X_{ij}$, converges to the wrong value $(1 - k^{-1})\sigma^2$. Although it is easy to correct the MLE in this particular situation, in general identifying the correction is a hard problem.

Ghosh proposed constructing an asymptotically efficient estimator for $\theta$ by viewing the $\xi_i$ as random variables arising from an unknown distribution $G$. The semiparametric model resulting from this can then be explored to find efficient estimators for $\theta$. In addition, efficiency in the Neyman–Scott model can be defined in terms of the semiparametric model so the two models have many interesting links between them. These links may be exploited and are studied in detail in [4], [5] and [6].

## 5. Bayesian inference

Ghosh has always been very fond of Bayesian ideas and, later in his career, he became more convinced that the Bayesian approach to statistics is more natural and fruitful. Over the course of his investigations, Ghosh examined each aspect of the Bayesian formulation, from construction of a prior to model selection, to asymptotic properties. And again, this can be seen as a continuation from his asymptotic work. After all, once the asymptotics are developed, we want to see how they can be used in a complete inference problem and the Bayesian setting provides a unified context. Indeed, Ghosh's contributions have helped speed the development of several branches of Bayesian analysis because of his asymptotic orientation.

Ghosh had always been pragmatic and thought that a good statistical method should have good frequentist properties as well as sensible conditional properties. Moreover, as in the frequentist case, asymptotics often play a vital role in Bayesian inference and one of the recurring themes in Ghosh's work has been the quest for frequentist properties of posterior distributions. As one of the leaders in developing objective Bayesian methods, he regularly worked to reconcile the two schools of thought. The paper [57] elaborately reviews issues and developments in objective Bayesian methodology.

Ghosh's Bayesian work can be broadly grouped into four categories. He has worked on frequentist matching and other objective priors. He has worked on determining the limiting behavior of posterior distributions in the parametric context. Then, he has turned his attention to richer model classes, examining Bayesian nonparametrics and model selection.

### *5.1. Matching and other objective priors*

Ghosh had never been very keen on the term noninformative to describe priors that are constructed through some automatic mechanism rather than through a subjective assessment of odds. His preference was to use these priors as objective or default priors in the absence of genuine subjective information. To him, such priors can be obtained by any one of various techniques including matching what a frequentist might use, invariance, entropy-type maximizations (reference priors), approximation, or anything that seemed reasonable.

The idea of a Bayesian choosing a prior so as to match frequentist inferences was originally introduced by Peers and Welch, but the term "probability matching prior" was first used by Ghosh and Mukerjee [42] and the approach became popular after Ghosh's presentation in the 1990 Valencia meeting. The basic idea is quite simple: choose a prior so that Bayesian notions like credibility approximately agree with the corresponding frequentist notions like confidence level. However, when asymptotic normality of the posterior distribution holds (discussed in the next subsection), it means that the variability according to the posterior distribution of a parameter is asymptotically equivalent to the sampling fluctuations of the MLE in the frequentist sense. This implies that Bayes-frequentist matching occurs for any prior under minimal restrictions. Consequently, to identify a prior uniquely, first order matching of limits is not enough. Satisfyingly, agreement continues to the next order, but only if the prior is of a certain form. Thus matching can be used to characterize a prior, which may then be thought of as objective at least in the sense that it was not chosen according to the personal views of the experimenter.

Of course, neither Bayesian credibility sets nor frequentist confidence sets are unique, so when $\theta$ is a scalar, it is natural to look at one sided intervals. If $W$ is a properly centered and normalized version of the parameter, then equating the posterior probability $P_\pi(W \leq t | X_1, \ldots, X_n)$ with the frequentist probability $P_\theta(W \leq t)$ for $t$ and ensuring both sides are $1 - \alpha$ up to $o(n^{-1/2})$ for each $\alpha$ leads to a differential equation. The Jeffreys prior is the solution to this equation.

A multiparameter version of this frequentist-Bayesian matching was used in [43]. In higher dimensions, the components of $W$ may be defined by successively computing the regression residual of the current component over the earlier components.

Naturally, this depends on the ordering of the parameters, but the dependence is not present when the parameters are orthogonal. In these cases, the matching criterion leads to partial differential equations. Curiously, Jeffreys' prior is a solution in some but not all cases: the location-scale problem is an important exception. In fact, it is well known that Jeffreys' prior, which is also the left Haar measure, may be an inappropriate choice in this case, so the matching criterion genuinely leads to sensible solutions even in high dimensional cases.

The matching prior is closely related to other important objective priors such as the reference prior. Reference priors often depend on the role of the parameter; nuisance parameters are treated differently from parameters of interest. Interestingly, in the two parameter case, a cute observation of Ghosh is that the reverse reference prior, rather than the reference prior itself, is probability matching. Here, by reverse reference prior, we mean the reference prior computed by reversing the roles of the parameter of interest and the nuisance parameter. More details and discussion of other properties, such as weak minimaxity, may be found in [42].

Although matching posterior probabilities does yield useful insight, highest posterior density (HPD) regions are more efficient credible sets from a Bayesian standpoint. Accordingly, matching the coverage probability of HPD regions with the credibility is an alternative that might be more appealing to some Bayesians. When this matching is done to $o(n^{-1})$, it leads to differential equations characterizing prior distributions. These were derived in [44]. In some cases, Jeffreys' prior solves these equations and so is a matching prior in the sense of coverage probability as well. A related paper is [46]. Matching the coverage of one-sided posterior credibility intervals for parametric functions up to $O(n^{-1})$ was studied in [17].

Alternatively, instead of characterizing a prior through matching, one might ask if there is some adjustment to make matching work for any prior satisfying mild general conditions. Indeed, in [47], it is shown, with examples, that if the center of the $(1-\alpha)$-HPD ellipsoid is appropriately shifted by a $o(n^{-1/2})$ amount, where the correction is obtained by solving an equation depending on the prior, then the resulting perturbed HPD ellipsoid's coverage is $1-\alpha + o(n^{-1})$.

Of course, there are many sensible notions of objectivity for a prior other than matching. Invariance is often the driving force in group models, where a group of transformation is acting on the parameter space and the parameter of interest is the maximal invariant parametric function. In [18], a detailed study of various priors such as the Chang–Eaves prior for group models is given in the light of matching and the marginalization paradox.

### 5.2. Limits of posterior distributions

One of the most intriguing results in statistics is the Bernstein–von Mises theorem, which states that the posterior distribution of the parameter centered at the MLE and scaled by $\sqrt{n}$ times the square root of the Fisher information converges to the standard normal distribution almost surely, as the sample size increases to infinity. This parallels the frequentist result that $\sqrt{n}(\hat{\theta} - \theta_{\text{true}})$ is asymptotically normal with variance given by the inverse Fisher information. In essence, posterior normality implies that in an asymptotic sense, at least to first order, any sensible Bayesian must agree eventually with frequentist notions of variability.

Ghosh worked to extend posterior normality in a variety of directions. One natural idea is to look at higher order properties so that the usual normal limit is viewed as merely the first term in an asymptotic expansion. This parallels the sense in which an Edgeworth expansion is an improvement over the standard central limit theorem. Johnson pioneered such expansions, but the probability statements in his expansions are in terms of the true distribution of the sample. Often, a Bayesian is more interested in bounds that are uniform on sets with high probability in the marginal distribution of the sample. In [61], precise conditions were given so that Johnson's expansion of the posterior distribution holds on a set with marginal probability $1 - O(n^{-r})$, where $r$ is the extra number of terms in the expansion, i.e., not counting the leading normal. It was also shown, by counterexamples, that some of the earlier published results in the field are incorrect.

Sometimes it is meaningful to condition on a statistic rather than the full data to obtain the posterior distribution. In particular, since the sample mean is a widely used summary measure, it is natural to ask if a version of the Bernstein–von Mises theorem holds when the posterior is computed given only the mean. Provided that expectation and variance are smooth functions, and the eigenvalues of the covariance matrix are uniformly bounded and bounded away from zero, it is shown in [15] that a normal limit for the posterior distribution is obtained. The variance of the limiting distribution can equal the variance of an observation, but in general, the normal limit can differ from that in the usual Bernstein–von Mises theorem, unless the sample mean is asymptotically sufficient. The proof is based on an Edgeworth expansion for the sample mean and a local limit theorem. The idea extends to independent but non-identically distributed observations.

More broadly, the Bernstein–von Mises phenomenon in a parametric family may be seen as the convergence of the posterior density of the standardized parameter to a non-degenerate distribution. In general, the centering need not be at the MLE, the scaling need not be $\sqrt{n}$ and the limit distribution need not be normal. Indeed, in some nonregular families such as the uniform distribution on $[0, \theta]$ or the location family of the exponential distribution, centering by the Bayes estimator and scaling by $n$ yields an exponential limit. This leads to the following question: For which families will a limit of the posterior distribution exist? When it does exist, what is the correct centering, scaling and limiting distribution? This problem is germane to approximating posterior distributions numerically when $n$ is large.

Under the general setting up of a parametric family considered by Ibragimov and Has'minskii in their book, a very elegant characterization was given in [35] and [23] in terms of the behavior of the limiting (local) likelihood ratio process of the model, $Z_n(u) = p(X^n; \theta + r_n u)/p(X^n; \theta)$, where $X^n$ is the observation at stage $n$ and $r_n$ is the appropriate normalizer for the problem. Usually $X^n = (X_1, \ldots, X_n)$ and $n$ is the sample size. Let $Z(u)$ stand for the weak limit of $Z_n(u)$ and $\xi(u) = Z(u)/\int Z(v)dv$, a random probability density. Under the natural scaling in the family, the posterior distribution converges to a limit, after appropriate centering, if and only if $\xi(u) = g(u + W)$ for some fixed probability density $g$ and a random variable $W$, i.e., as a random element in $L_1$, $\xi(\cdot)$ is a random location shift of a fixed probability density $g$. When this holds, $g$ is the limit of the posterior density. Clearly, this is a stringent representation, so in many nonregular cases the posterior distribution will not have a limit.

Interestingly, in the regular cases, local asymptotic normality implies that $\xi(u) = g(u + W)$, in which $g$ is normal and $W$ is a random normal shift. Thus this yields a Bernstein–von Mises theorem under an extremely general condition. A similar limit theorem holds with an exponential limit whenever densities are positively supported on an interval $[a(\theta), b(\theta)]$, where the support is either expanding or contracting.

While it is disappointing to find that posterior limits exist only in relatively rare cases, it does not rule out the possibility of finding useful approximation to the posterior distributions depending on the sample size $n$. For change-point problems, where the density jumps from a positive value to another positive value at an unknown location but is otherwise smooth, a useful approximation was obtained in [24] by normalizing an approximation to the likelihood ratio process. It turns out that a certain mixture of $n$ many truncated and shifted exponential densities is a good approximation.

### *5.3. Bayesian nonparametrics*

Ghosh's involvement with Bayesian nonparametrics started in the mid 90's with the paper [51] attempting to determine whether the priors used for survival analysis lead to consistency under censoring. This paper showed that for the Dirichlet process, the posterior under censoring can be represented as a Pòlya tree process whose partition depends on the data, and then consistency can be obtained from the tail-free property of Pòlya tree processes. The question is followed up in subsequent papers [53] and [36]. Since then, Ghosh has continued to be one of the most important contributors to understanding the asymptotics of Bayesian nonparametrics.

For instance, the search for a noninformative prior for infinite dimensional models, as an extension to the finite dimensional case, is ongoing. One approach is to generalize the notion of a uniform distribution. This was proposed in [19] using uniform distributions on discrete approximations to a space found by maximal $\epsilon$-dispersed sets. Even in the parametric setting this approach is fundamental and leads to Jeffreys' prior. The approach gives consistency in infinite-dimensional cases.

More typically, Ghosh was strongly motivated by the examples of inconsistency of posterior distributions in infinite dimensional models. While he appreciated those illuminating examples, he was always hopeful that Bayes' methods would work if priors were constructed properly. He was particularly fond of the Kullback–Leibler property which requires that the true distribution be in the support of the prior when distances are measured in terms of Kullback–Leibler divergence. That is, the prior should assign strictly positive probability to every Kullback–Leibler neighborhood around the true distribution.

Because of this, Ghosh thought the Dirichlet process was inappropriate in many contexts, despite its evident utility, since its discreteness means it fails to have anything in its Kullback–Leibler support. In [20], it was shown that a prior with the Kullback–Leibler property, such as a suitable Pòlya tree or a Dirichlet mixture process, can overcome the inconsistency property of Dirichlet processes for estimating a location parameter. Essentially the same phenomenon appears in linear regression models as shown in [1]. In that paper, the first extension of a general posterior consistency theorem to independent non-identically distributed variables is also developed.

In Bayesian nonparametrics, consistency often combines testing concepts with sieves. A celebrated result of Schwartz emphasizes the role of tests for the true density $f_0$ versus the complement of a neighborhood, say $V$, around it. The basic idea is to construct tests, by covering $V^c$ with many small balls, say $B_i$'s, and testing $f_0$ versus $B_i$ for each $i$ using powerful tests. One can then simply look at the maximum of all tests against each small ball, whose type II error probability is clearly under control and the type I error probability bounded by the common exponential bound for error probability multiplied by the number of small balls required to cover $V^c$. Thus the concept of metric entropy, which is the logarithm of the number of balls required to cover a set, comes into the picture. Generally $V^c$ is not compact, and it is not possible to cover it by finitely many small balls. The difficulty can be overcome by using a sieve, which is a sequence of increasing subsets of a parameter space that gradually fill out the whole parameter space. One may ignore the portion of the parameter space outside the sieve as long as that part has exponentially small prior probability. Now one must control the metric entropy of the sieve to ensure that it does not grow faster than a small multiple of $n$. This style of proof gives consistency for density estimation with Dirichlet mixtures of normal kernels as shown in [21], providing a large sample justification for the most widely used Bayesian density estimator.

This approach works for density estimation with other priors in place of the Dirichlet mixtures. In [68], consistency is obtained for the logistic Gaussian prior for a density, that is, a prior on densities obtained by exponentiating and then normalizing a Gaussian process.

The importance of entropy for posterior consistency appeared in [21]. There it is seen that in the nonparametric setting prior positivity at the true density must be satisfied, but in terms of special neighborhoods given by the Kullback–Leibler number. Moreover, it must be possible to choose a sieve whose entropy grows no faster than the rate $O(n)$, while ensuring that the prior probability of the complement of sieve is exponentially small as $n$ increases. This observation led to the derivation of the results on posterior convergence rates in [25] in the sense that the conditions for rates can be viewed as quantitative analogs of the conditions for consistency.

For instance, instead of just requiring that the prior for a fixed $\epsilon$ neighborhood in the Kullback–Leibler sense has positive probability, one now needs to show that the prior probability of the Kullback–Leibler neighborhood of radius $\epsilon_n$ is at least $e^{-n\epsilon_n^2}$, where $\epsilon_n$ is the intended rate of posterior convergence. In a similar manner, requiring that the $\epsilon_n$- entropy of the sieve be bounded by a multiple of $n\epsilon_n^2$ is also reminiscent of the condition that the $\epsilon$-entropy of the sieve should be bounded by a small multiple of $n$. Thus, for fixed $\epsilon_n$, this reduces to the condition for consistency.

The paper also constructs a prior achieving optimal rates of convergence by bracketing densities above and below by two functions – choosing a finite collection of brackets to provide upper and lower bounds for any probability density in the given class. This ensures good approximation of any function within the bracket together with a control over likelihood ratios. This can be viewed as a refinement of the construction proposed in [19]. Other approaches to optimal rates are also discussed, most notably, through exponential families generated using a B-spline basis.

Many aspects of Bayesian asymptotics for infinite dimensional models are neatly summarized in the review [22], and thoroughly discussed in [52], which to date is the only book dealing with asymptotic results in Bayesian nonparametrics.

### *5.4. Model selection and Bayesian hypothesis testing*

Testing hypotheses is a major area where frequentist and Bayesian procedures often differ substantially. There is a tendency for frequentist methods to over-reject just as there is a tendency for Bayes' methods to under-reject, as in the Lindley paradox. Results such as the consistency of the Bayesian information criterion (BIC) bridge the gap somewhat because the BIC approximates Bayes factors and is frequentist consistent for model selection under appropriate conditions in the sense that the BIC selects the correct model with probability tending to one.

These properties of the BIC are valid only if the dimension $p$ of the model remains fixed. However, for many applications, especially for complex data containing numerous variables commonly arising nowadays, the BIC fails to approximate the Bayes factor adequately and is consistent. The main reason for the failure is ignoring certain terms in an expansion of the Bayes factor which are not negligible when $p \to \infty$. The difficulty can be avoided by paying proper attention to these terms. In [3], a correction is proposed by introducing two more terms, one is proportional to $p$ and the other to $\log p$, as well as changing the meaning of the sample size to the number of replications. The resulting "generalized BIC" then selects the correct model with increasingly high probability. Another generalization of BIC is developed in [11] which works in a general exponential family. These generalizations of BIC are powerful tools to overcome the challenges posed by high-dimensional data problems of contemporary statistics.

In model selection problems, the definition of optimality is often tricky. An appealing approach is comparison with the oracle, that is, with the best procedure (for a given loss function) which uses the knowledge of the correct model in making decisions. A parametric empirical Bayes (PEB) approach approximates the Bayes factor by deriving the rule in a parametric model but estimating the parameters in the penalty function by a penalized likelihood criterion with data dependent penalty. In [66], the relative performance of a PEB, the AIC and the BIC were thoroughly studied through asymptotics and simulations under both 0-1 and prediction loss. The conclusion is that the BIC performs badly, but PEB rules perform quite satisfactorily, and so does the AIC. If Bayes estimates are used in making predictions, instead of least squares estimators, a PEB performs better than the AIC.

One particular difficulty with the Bayes factor is that it is undefined when improper priors are used in individual models. Various remedies are proposed in the literature, based on the idea of using a part of the information contained in the data (training portion) to make priors proper and use the remaining portion in Bayesian analysis with the obtained "proper prior". Since this typically depends on the ordering of the data, some kind of averaging, through bootstrap or cross validation, over different choices of the training portion is desired. A particularly popular candidate among these Bayes factors is obtained by taking a geometric average. In [67], such Bayes factors are studied through asymptotics as the proportion of the training sample varies, and conditions for consistency are obtained as the total sample size goes to infinity. It turns of that predictive optimality of the "geometric Bayes factor" as it is often claimed is not entirely correct.

There are many other significant papers on model selection authored by Ghosh. In [10], optimality of the AIC in inference about Brownian motion is shown. The reviews [56] and [55] contain wealth of information on Bayesian model selection.

## 6. Concluding remarks

Overall, Ghosh's work in statistics reveals a progression. He began with individual data points, proceeded to data summarization, and then to the asymptotics of inference. Ghosh's results there were a successful attempt to map out where the accumulation of data tend to point. In a sense, asymptotic limits are the ultimate data summarization. Then, putting it all together, Ghosh turned to the Bayesian formulation, examining each of its components, prior, model, posterior, in turn, to permit a comprehensive and unified study of the statistical problem. Indeed, his recent work on Bayesian nonparametrics is a further generalization, again a logical step because it builds on his earlier work by using ever richer model classes.

In fact, Ghosh has worked in many more areas of statistics, apart from those outlined above, as well as working on a variety of applications. These topics include distribution theory, decision theory, robustness, finite population sampling, reliability, quality control, modeling hydrocarbon discovery, geological mapping and DNA fingerprinting.

Finally, every great researcher has a strategy, a method or a drive, often summarized in a maxim, that guides or motivates their intellectual endeavors. One of Ghosh's maxims was the injunction: "Settle the question!" By this he meant formulate a question so that answering it gives you something definite for the formulation of another question. As can be inferred from the progression of his work, his questioning led him to an ever broader view of the statistical problem, culminating in a Bayesian treatment of high-dimensional models, nonparametric or not. Ghosh's injunction to settle questions has helped, and will continue to help, researchers all over the world to think deeply about the most important issues.

## References

[1] AMEWOU-ATISSO, M., GHOSAL, S., GHOSH, J. K. AND RAMAMOORTHI, R. V. (2003). Posterior consistency for semi-parametric regression problems. *Bernoulli* **9** 291–312. MR1997031

[2] BASU, D. AND GHOSH, J. K. (1969). Sufficient statistics in sampling from a finite universe. *Bull. Inst. Internat. Statist.* **42** 850–858. MR0286197

[3] BERGER, J. O., GHOSH, J. K. AND MUKHOPADHYAY, N. (2003). Approximations and consistency of Bayes factors as model dimension grows. *J. Statist. Plann. Inference* **112** 241–258. MR1961733

[4] BHANJA, J. AND GHOSH, J. K. (1992). Efficient estimation with many nuisance parameters. I. *Sankhyā Ser. A* **54** 1–39. MR1189781

[5] BHANJA, J. AND GHOSH, J. K. (1992). Efficient estimation with many nuisance parameters. II. *Sankhyā Ser. A* **54** 135–156. MR1192091

[6] BHANJA, J. AND GHOSH, J. K. (1992). Efficient estimation with many nuisance parameters. III. *Sankhyā Ser. A* **54** 297–308. MR1216288

[7] BHATTACHARYA, R. N. AND GHOSH, J. K. (1978). On the validity of the formal Edgeworth expansion. *Ann. Statist.* **6** 434–451. MR0471142

[8] BHATTACHARYA, R. N. AND GHOSH, J. K. (1988). On moment conditions for valid formal Edgeworth expansions. *J. Multivariate Anal.* **27** 68–79. MR0971173

[9] BICKEL, P. J. AND GHOSH, J. K. (1990). A decomposition for the likelihood ratio statistic and the Bartlett correction – a Bayesian argument. *Ann. Statist.* **18** 1070–1090. MR1062699

[10] CHAKRABARTI, A. AND GHOSH, J. K. (2006). Optimality of AIC in inference about Brownian motion. *Ann. Inst. Statist. Math.* **58** 1–20. MR2281204

[11] CHAKRABARTI, A. AND GHOSH, J. K. (2006). A generalization of BIC for the general exponential family. *J. Statist. Plann. Inference* **136** 2847–2872. MR2281234

[12] CHANDRA, T. K. AND GHOSH, J. K. (1978). Comparison of tests with same Bahadur efficiency. *Sankhyā Ser. A* **40** 253–277. MR0589281

[13] CHANDRA, T. K. AND GHOSH, J. K. (1979). Valid asymptotic expansions for the likelihood ratio statistic and other perturbed chi-square variables. *Sankhyā Ser. A* **41** 22–47. MR0615038

[14] Chandra, T. K. and Ghosh, J. K. (1980). Valid asymptotic expansions for the likelihood ratio and other statistics under contiguous alternatives. *Sankhyā Ser. A* **42** 170–184. MR0656254

[15] Clarke, B. and Ghosh, J. K. (1995). Posterior convergence given the mean. *Ann. Statist.* **23** 2116–2144. MR1389868

[16] DasGupta, A. and Ghosh, J. K. (1983). Some remarks on second-order admissibility in the multiparameter case. *Sankhyā Ser. A* **45** 181–190. MR0748457

[17] Datta, G. S. and Ghosh, J. K. (1995). On priors providing frequentist validity for Bayesian inference. *Biometrika* **82** 37–45. MR1332838

[18] Datta, G. S. and Ghosh, J. K. (1995). Noninformative priors for maximal invariant parameter in group models. *Test* **4** 95–114. MR1365042

[19] Ghosal, S., Ghosh, J. K. and Ramamoorthi, R. V. (1997). Noninformative priors via sieves and packing numbers. In *Advances in Statistical Decision Theory and Applications* 119–132. *Stat. Ind. Technol.* Birkhäuser, Boston. MR1479180

[20] Ghosal, S., Ghosh, J. K. and Ramamoorthi, R. V. (1999). Consistent semiparametric Bayesian inference about a location parameter. *J. Statist. Plann. Inference* **77** 181–193. MR1687955

[21] Ghosal, S., Ghosh, J. K. and Ramamoorthi, R. V. (1999). Posterior consistency of Dirichlet mixtures in density estimation. *Ann. Statist.* **27** 143–158. MR1701105

[22] Ghosal, S., Ghosh, J. K. and Ramamoorthi, R. V. (1999). Consistency issues in Bayesian nonparametrics. In *Asymptotics, Nonparametrics, and Time Series* 639–667. *Statist. Textbooks Monogr.* **158**. Dekker, New York. MR1724711

[23] Ghosal, S., Ghosh, J. K. and Samanta, T. (1995). On convergence of posterior distributions. *Ann. Statist.* **23** 2145–2152. MR1389869

[24] Ghosal, S., Ghosh, J. K. and Samanta, T. (1999). Approximation of the posterior distribution in a change-point problem. *Ann. Inst. Statist. Math.* **51** 479–497. MR1722841

[25] Ghosal, S., Ghosh, J. K. and van der Vaart, A. W. (2000). Convergence rates of posterior distributions. *Ann. Statist.* **28** 500–531. MR1790007

[26] Ghosh, J. K. (1960). On some properties of sequential *t*-test. *Calcutta Statist. Assoc. Bull.* **9** 77–86. MR0114277

[27] Ghosh, J. K. (1960). On the monotonicity of the *OC* of a class of sequential probability ratio tests. *Calcutta Statist. Assoc. Bull.* **9** 139–144. MR0117849

[28] Ghosh, J. K. (1961). On the optimality of probability ratio tests in sequential and multiple sampling. *Calcutta Statist. Assoc. Bull.* **10** 73–92. MR0130774

[29] Ghosh, J. K. (1964). Bayes solutions in sequential problems for two or more terminal decisions and related results. *Calcutta Statist. Assoc. Bull.* **13** 101–122. MR0172422

[30] Ghosh, J. K. (1971). A new proof of the Bahadur representation of quantiles and an application. *Ann. Math. Statist.* **42** 1957–1961. MR0297071

[31] Ghosh, J. K. (1991). Higher order asymptotics for the likelihood ratio, Rao's and Wald's tests. *Statist. Probab. Lett.* **12** 505–509. MR1143747

[32] Ghosh, J. K. (1994). *Higher Order Asymptotics. NSF-CBMS Regional Conference in Probability and Statistics* **4**. IMS, Hayward, CA.

[33] Ghosh, J. K. and Chaudhuri, A. R. (1984). An invariant SPRT for identification. *Sequential Anal.* **3** 99–120. MR0767249

[34] Ghosh, J. K., Delampady, M. and Samanta, T. (2007). *An Introduction*

*to Bayesian Analysis, Theory and Methods.* Springer, New York. MR2247439
[35] GHOSH, J. K., GHOSAL, S. AND SAMANTA, T. (1994). Stability and convergence of the posterior in non-regular problems. In *Statistical Decision Theory and Related Topics* **V** 183–199. Springer, New York. MR1286304
[36] GHOSH, J. K., HJORT, N. L., MESSAN, C. AND RAMAMOORTHI, R. V. (2006). Bayesian bivariate survival estimation. *J. Statist. Plann. Inference* **136** 2297–2308. MR2235060
[37] GHOSH, J. K., MORIMOTO, H. AND YAMADA, S. (1981). Neyman factorization and minimality of pairwise sufficient subfields. *Ann. Statist.* **9** 514–530. MR0615428
[38] GHOSH, J. K. AND MUKERJEE, R. (1989). Some optimality results on Stein's two-stage sampling. In *Statistical Data Analysis and Inference* 251–256. North-Holland, Amsterdam. MR1089640
[39] GHOSH, J. K. AND MUKERJEE, R. (1990). Improvement in Stein's procedure using a random confidence coefficient. *Calcutta Statist. Assoc. Bull.* **40** 145–152. MR1172640
[40] GHOSH, J. K. AND MUKERJEE, R. (1991). Characterization of priors under which Bayesian and frequentist Bartlett corrections are equivalent in the multiparameter case. *J. Multivariate Anal.* **38** 385–393. MR1131727
[41] GHOSH, J. K. AND MUKERJEE, R. (1992). Bayesian and frequentist Bartlett corrections for likelihood ratio and conditional likelihood ratio tests. *J. Roy. Statist. Soc. Ser. B* **54** 867–875. MR1185228
[42] GHOSH, J. K. AND MUKERJEE, R. (1992). Non-informative priors. In *Bayesian Statistics* **4** 195–210. Oxford Univ. Press, New York. MR1380277
[43] GHOSH, J. K. AND MUKERJEE, R. (1993). On priors that match posterior and frequentist distribution functions. *Canad. J. Statist.* **21** 89–96. MR1221860
[44] GHOSH, J. K. AND MUKERJEE, R. (1993). Frequentist validity of highest posterior density regions in the multiparameter case. *Ann. Inst. Statist. Math.* **45** 293–302. MR1232496
[45] GHOSH, J. K. AND MUKERJEE, R. (1994). Adjusted versus conditional likelihood: power properties and Bartlett-type adjustment. *J. Roy. Statist. Soc. Ser. B* **56** 185–188. MR1257806
[46] GHOSH, J. K. AND MUKERJEE, R. (1995). Frequentist validity of highest posterior density regions in the presence of nuisance parameters. *Statist. Decisions* **13** 131–139. MR1342734
[47] GHOSH, J. K. AND MUKERJEE, R. (1995). On perturbed ellipsoidal and highest posterior density regions with approximate frequentist validity. *J. Roy. Statist. Soc. Ser. B* **57** 761–769. MR1354080
[48] GHOSH, J. K. AND MUKERJEE, R. (2001). Test statistics arising from quasi likelihood: Bartlett adjustment and higher-order power. *J. Statist. Plann. Inference* **97** 45–55. MR1851373
[49] GHOSH, J. K., MUKERJEE, R. AND SEN, P. K. (1996). Second-order Pitman admissibility and Pitman closeness: the multiparameter case and Stein-rule estimators. *J. Multivariate Anal.* **57** 52–68. MR1392577
[50] GHOSH, J. K., PURKAYASTHA, S. AND SAMANTA, T. (2004). Sequential probability ratio tests based on improper priors. *Sequential Anal.* **23** 585–602. MR2103910
[51] GHOSH, J. K. AND RAMAMOORTHI, R. V. (1995). Consistency of Bayesian inference for survival analysis with or without censoring. In *Analysis of Censored Data. IMS Lecture Notes Monogr. Ser.* **27**. IMS, Hayward, CA. MR1483342
[52] GHOSH, J. K. AND RAMAMOORTHI, R. V. (2003). *Bayesian Nonparametrics.*

Springer, New York. MR1992245
- [53] GHOSH, J. K., RAMAMOORTHI, R. V. AND SRIKANTH, K. R. (1999). Bayesian analysis of censored data. *Statist. Probab. Lett.* **41** 255–265. MR1672393
- [54] GHOSH, J. K. AND ROY, K. K. (1972). Families of densities with non-constant carriers which have finite dimensional sufficient statistics. *Sankhyā Ser. A* **34** 205–226. MR0378158
- [55] GHOSH, J. K. AND SAMANTA, T. (2001). Model selection – an overview. *Current Science* **80** (9), 1135–1144.
- [56] GHOSH, J. K. AND SAMANTA, T. (2002). Nonsubjective Bayes testing – an overview. *J. Statist. Plann. Inference* **103** 205–223. MR1896993
- [57] GHOSH, J. K. AND SAMANTA, T. (2002). Towards a nonsubjective Bayesian paradigm. *Uncertainty and Optimality* 1–69. World Sci. Publ., River Edge, NJ. MR1955963
- [58] GHOSH, J. K., SEN, P. K. AND MUKERJEE, R. (1994). Second-order Pitman closeness and Pitman admissibility. *Ann. Statist.* **22** 1133–1141. MR1311968
- [59] GHOSH, J. K. AND SINHA, B. K. (1981). A necessary and sufficient condition for second order admissibility with applications to Berkson's bioassay problem. *Ann. Statist.* **9** 1334–1338. MR0630116
- [60] GHOSH, J. K. AND SINHA, B. K. (1982). Third order efficiency of the MLE – a counterexample. *Calcutta Statist. Assoc. Bull.* **31** 151–158. MR0702402
- [61] GHOSH, J. K., SINHA, B. K. AND JOSHI, S. N. (1982). Expansions for posterior probability and integrated Bayes risk. In *Statistical Decision Theory and Related Topics III* **1** 403–456. Academic Press, New York. MR0705299
- [62] GHOSH, J. K., SINHA, B. K. AND SUBRAMANYAM, K. (1979). Edgeworth expansions for Fisher-consistent estimators and second order efficiency. *Calcutta Statist. Assoc. Bull.* **28** 1–18. MR0586079
- [63] GHOSH, J. K., SINHA, B. K. AND WIEAND, H. S. (1980). Second order efficiency of the MLE with respect to any bounded bowl-shaped loss function. *Ann. Statist.* **8** 506–521. MR0568717
- [64] GHOSH, J. K. AND SUBRAMANYAM, K. (1974). Second order efficiency of maximum likelihood estimators. *Sankhyā Ser. A* **36** 325–358. MR0428572
- [65] HALL, W. J., WIJSMAN, R. A. AND GHOSH, J. K. (1965). The relationship between sufficiency and invariance with applications in sequential analysis. *Ann. Math. Statist.* **36** 575–614. MR0178552
- [66] MUKHOPADHYAY, N. AND GHOSH, J. K. (2003). Parametric empirical Bayes model selection – some theory, methods and simulation. In *Probability, Statistics and Their Applications: Papers in Honor of Rabi Bhattacharya* 229–245. *IMS Lecture Notes Monogr. Ser.* **41**. IMS, Beachwood, OH. MR1999424
- [67] MUKHOPADHYAY, N., GHOSH, J. K. AND BERGER, J. O. (2005). Some Bayesian predictive approaches to model selection. *Statist. Probab. Lett.* **73** 369–379. MR2187852
- [68] TOKDAR, S. T. AND GHOSH, J. K. (2007). Posterior consistency of logistic Gaussian process priors in density estimation. *J. Statist. Plann. Inference* **137** 34–42. MR2292838